\documentclass[%
 reprint,
 amsmath,amssymb,
 aps, 
pre,
]{revtex4-1}

\usepackage{float}
\usepackage{subfigure}

\usepackage{graphicx}
\usepackage{dcolumn}
\usepackage{bm}

\begin{document}

\preprint{APS/123-QED}

\title{Nonlinear Dynamics of the Rock-Paper-Scissors Game with Mutations}

\author{Danielle F. P. Toupo}
\author{Steven H. Strogatz}%
\affiliation{%
Center for Applied Mathematics, Cornell University, Ithaca, New York 14853\\
}%

\date{\today}

\begin{abstract}
We analyze the replicator-mutator equations for the Rock-Paper-Scissors game. Various graph-theoretic patterns of mutation are considered, ranging from a single unidirectional mutation pathway between two of the species, to global bidirectional mutation among all the species. Our main result is that the coexistence state, in which all three species exist in equilibrium, can be destabilized by arbitrarily small mutation rates. After it  loses stability, the coexistence state gives birth to a stable limit cycle solution created in a supercritical Hopf bifurcation. This attracting periodic solution exists for all the mutation patterns considered, and persists arbitrarily close to the limit of zero mutation rate and a zero-sum game.  

\end{abstract}

\maketitle


\section{Introduction}
In the children's game of Rock-Paper-Scissors, paper wraps rock, rock smashes scissors, and scissors cut paper. Theoretical biologists  have used this game as a metaphor in ecology and evolutionary biology to describe interactions among three competing species in which each species has an advantage over one of its opponents but not the other~\cite{hofbauer1998evolutionary,nowak2006evolutionary,szabo2007evolutionary,szolnoki2014cyclic}, 
a situation that can result in cyclic dominance~\cite{hofbauer1998evolutionary,szolnoki2014cyclic,kerr2002local,kirkup2004antibiotic,mobilia2010oscillatory,reichenbach2007mobility,szabo2004rock,szolnoki2009phase,szczesny2013does}. Interactions of this type occur between three competing strains of \textit{E. coli}~\cite{kerr2002local} as well as in the mating strategies of side-blotched lizards~\cite{sinervo1996rock}.  The Rock-Paper-Scissors game also has applications outside of biology. For example, it has been used to analyze the dynamics of a sociological system with three strategies in a public goods game~\cite{semmann2003volunteering}. 

Mathematically, the dynamics of the Rock-Paper-Scissors game is often studied using the replicator equations~\cite{hofbauer1998evolutionary,nowak2006evolutionary}, a system of coupled nonlinear differential equations that have been applied in diverse scientific settings. For example, replicator equations have been used to model the evolution of language, fashion, autocatalytic chemical networks, behavioral dynamics, and multi-agent decision making in social networks~\cite{hofbauer1998evolutionary,nowak2006evolutionary,szabo2007evolutionary,szolnoki2014cyclic,galla2011cycles,imhof2005evolutionary,mitchener2004chaos,pais2012Hopf,pais2011limit,stadler1992mutation,szabo2004cooperation}.

In the specific context of the Rock-Paper-Scissors game, the replicator equations are most often studied in the absence of mutation~\cite{hofbauer1998evolutionary,may1975nonlinear,nowak2006evolutionary,szabo2007evolutionary}. The dynamics in that case tend to exhibit one of three types of long-term behavior, depending on a parameter $\epsilon$ that characterizes how far the game is from being a zero-sum game. The three types of behavior are: (i) stable coexistence of all three species, (ii) neutrally stable cycles that fill the whole state space, and (iii) large-amplitude heteroclinic cycles in which each species in turn almost takes over the whole population and then almost goes extinct.  

In 2010, Mobilia~\cite{mobilia2010oscillatory} broke new ground by asking what would happen if the strategies are allowed to mutate into one another. For simplicity he assumed that every species could mutate into every other at a uniform rate $\mu$. For this highly symmetrical case (which we call global mutation), Mobilia~\cite{mobilia2010oscillatory} found a new form of behavior, not present in the mutation-free case:  the system could settle into a stable limit-cycle oscillation, born from a supercritical Hopf bifurcation.

In this paper we investigate what happens in other mutational regimes. For example, we consider the dynamics when precisely one species mutates into one other, or when two species mutate into two others, as well as more complicated patterns of mutation. In every case we find that stable limit cycles can occur, and we calculate the regions in parameter space where such attracting periodic behavior occurs.  

Our work was motivated by a question that arose in one of our previous studies~\cite{toupo}. There we had explored the dynamics of the repeated Prisoner's Dilemma game for three strategies: always defect (ALLD), always cooperate (ALLC), and tit-for-tat (TFT).
We incorporated mutations at a rate $\mu$   and a complexity cost $c$ of playing TFT into the replicator equations and analyzed the six possible single-mutation cases (exactly one strategy mutates into one other with rate $\mu$) as well as the global mutation case (each strategy mutates into the two others at the same rate $\mu$). Our results showed that stable limit cycles of cooperation and defection were possible for \emph{every} pattern of mutation we considered. What was particularly striking was that stable cycles occurred for parameter values arbitrarily close to the structurally unstable limiting case of zero mutation rate and zero complexity cost $(\mu = 0, c=0)$. We conjectured, but were unable to prove, that stable limit cycles would continue to exist near this point for \emph{any} pattern of mutation, not just the single-mutation patterns we looked at in Ref.~\cite{toupo}. 

In hopes of finding a more tractable system where the same phenomena might occur, we turned to the Rock-Paper-Scissors game. The analogous question is, does this system always exhibit stable limit cycles for parameters arbitrarily close to the zero-sum, zero-mutation-rate limit, for any possible pattern of mutation? For single mutations, the answer is yes, as we will show below. For arbitrary patterns of mutation, the answer again appears to be yes, but we have only managed to prove this under a further constraint, namely that the mutation pattern preserves the symmetrical coexistence state where all three strategies are equally populated.


\section{Rock-Paper-Scissors game}
The standard Rock-Paper-Scissors game is a zero-sum game with payoff matrix given in Table~\ref{tab:payoff0sum}. The entries show the payoff received by the row player when playing with the column player. 

\begin{table}[h!]
\caption{\label{tab:payoff0sum}%
Payoff matrix of a zero-sum Rock-Paper-Scissors game
}
\begin{ruledtabular}
\begin{tabular}{lcdr}
\textrm{}&
\textrm{Rock}&
\multicolumn{1}{c}{Paper}&
\textrm{Scissors}\\
\colrule
Rock & 0 & -1 & 1\\
Paper & 1 & 0 & -1\\
Scissors & -1 & 1 & 0\\
\end{tabular}
\end{ruledtabular}
\end{table}

According to the payoff matrix in Table~\ref{tab:payoff0sum}, each species gets a payoff  $0$ when playing against itself. When playing against a different species, the winner gets a payoff $1$ while the loser gets $-1$. 

A more general payoff matrix, considered in Ref.~\cite{mobilia2010oscillatory}, allows for non-zero-sum games. Now the winner gets a payoff $1$ and the loser gets $-\epsilon$, as shown in Table~\ref{tab:payoffeps1}.  This  game is zero sum if and only if $\epsilon = 1$. 

\begin{table}[h!]
\caption{\label{tab:payoffeps1}%
Payoff matrix of a Rock-Paper-Scissors game with parameter $\epsilon$. The game is a zero-sum game when $\epsilon=1$.
}
\begin{ruledtabular}
\begin{tabular}{lcdr}
\textrm{}&
\textrm{Rock}&
\multicolumn{1}{c}{Paper}&
\textrm{Scissors}\\
\colrule
Rock & 0 &-\epsilon & 1\\
Paper & 1 & 0 & $-\epsilon$\\
Scissors & $-\epsilon$ & 1 & 0\\
\end{tabular}
\end{ruledtabular}
\end{table}

For convenience, we redefine the entries of the payoff matrix so that the zero-sum case corresponds to $\epsilon =0$ rather than $\epsilon = 1$. The payoff matrix for the rest of this paper is shown in Table~\ref{tab:payoffeps0}. 

\begin{table}[h!]
\caption{\label{tab:payoffeps0}%
Payoff matrix of a Rock-Paper-Scissors game, redefined so that the zero-sum game has $\epsilon=0.$
}
\begin{ruledtabular}
\begin{tabular}{lcdr}
\textrm{}&
\textrm{Rock}&
\multicolumn{1}{c}{Paper}&
\textrm{Scissors}\\
\colrule
Rock & 0 &-(\epsilon+1) & 1\\
Paper & 1 & 0 & $-(\epsilon+1)$\\
Scissors & $-(\epsilon+1)$ & 1 & 0\\
\end{tabular}
\end{ruledtabular}
\end{table}


\section{Global mutations}
\begin{figure}[h!]
\hspace{0.5cm}\includegraphics[trim = 0mm 22cm 0mm 2cm,clip,scale = 0.4]{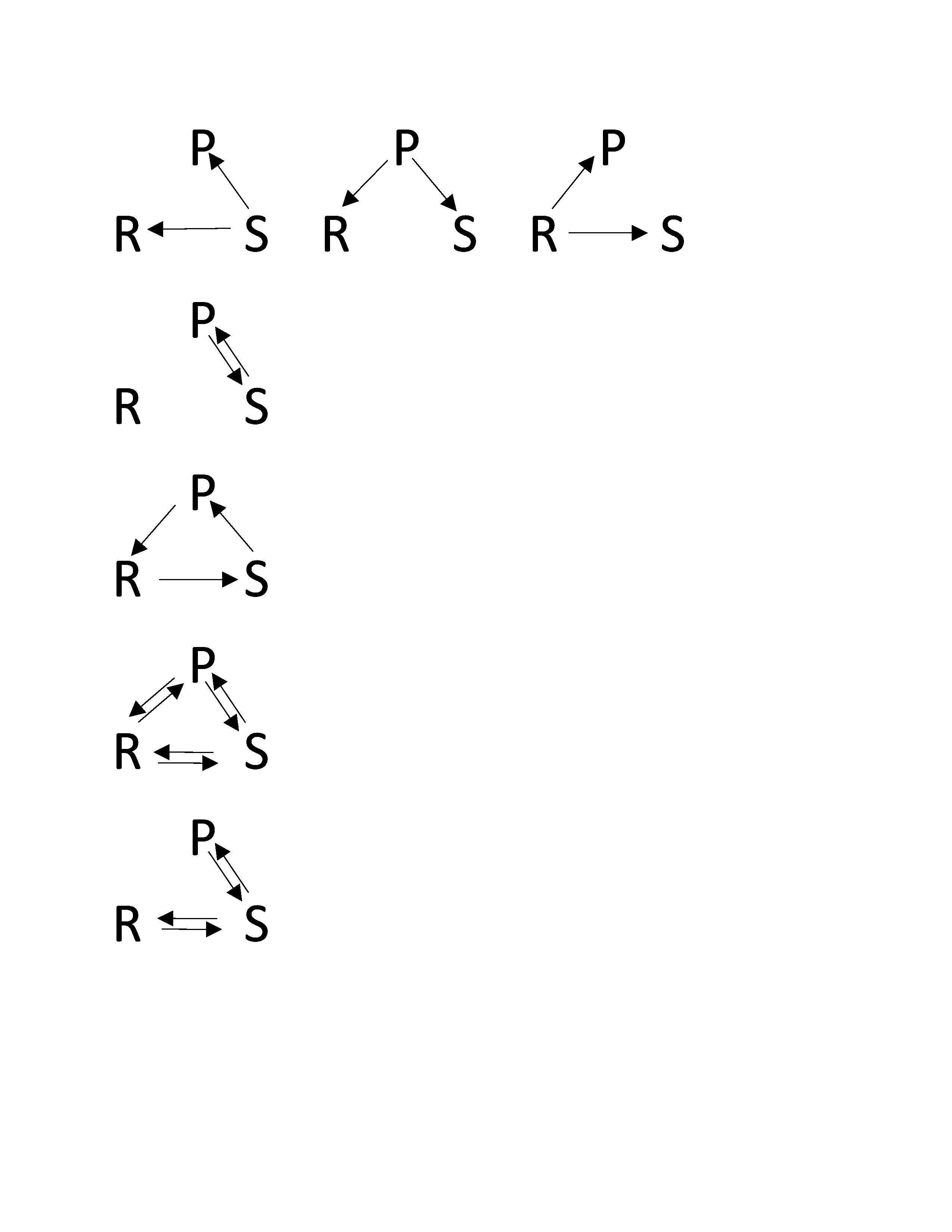}
\caption{\label{fig:arrows} Pathways for global mutation.}
\end{figure}

The most symmetrical Rock-Paper-Scissors game with mutation is the one with global mutation, where each species mutates into the other two with a rate $\mu$, as shown in Fig.~\ref{fig:arrows}.


Suppose that global mutation occurs in a well-mixed population with $N$ individuals. Let $x, y$ and $z$ denote the relative frequencies of individuals playing rock, paper and scissors, respectively. Then $x+y+z=1$ or $z=1-x-y$. By eliminating $z$ in this fashion, one can capture the dynamics of the three strategies by studying $x$ and $y$ alone. Following Ref.~\cite{mobilia2010oscillatory}, the replicator-mutator equations can be written as 

\begin{eqnarray}
\dot{x} &=& x \left(f_x - \phi \right) +\mu \left(-2 x+y+z\right)  \nonumber\\
\dot{y} &=& y   \left(f_y - \phi \right) +\mu \left(-2y+x+z \right)   \label{eqn:globalrighttriangle}
\end{eqnarray}

\noindent where $f_i$ is the fitness of strategy $i$, defined as its expected payoff against the current mix of strategies,
and $\phi = xf_x+yf_y+zf_z$ is the average fitness in the whole population. For the payoff matrix defined above in Table~\ref{tab:payoffeps0}, and using $z=1-x-y$, we find $$f_x = 1- x - (\epsilon +2) y $$ and $$f_y = ( \epsilon + 2) x + (\epsilon + 1) (y - 1).$$

A conceptual disadvantage of eliminating $z$ in favor of $x$ and $y$ is that the system's cyclic symmetry becomes less apparent. To highlight it, one can plot the phase portraits of the system on the equilateral triangle defined by the face of the simplex $x+y+z=1$, where $0 \leq x, y, z \leq 1$. This mapping of the phase portrait onto the equilateral triangle can be achieved by using the following transformation:
\begin{equation}
 \binom{X}{Y} =  \begin{pmatrix} 1 & \frac{1}{2}  \\ 0 & \frac{\sqrt{3}}{2} \end{pmatrix} \binom{x}{y} \label{eqn:transformation}.
\end{equation}
In what follows, all phase portraits will be plotted in $(X,Y)$ space, but we will still indicate the values of $x,y,z$ at the vertices of the simplex, since these variables have clearer interpretations.

\begin{figure}[h!]

\center
\includegraphics[scale=0.5]{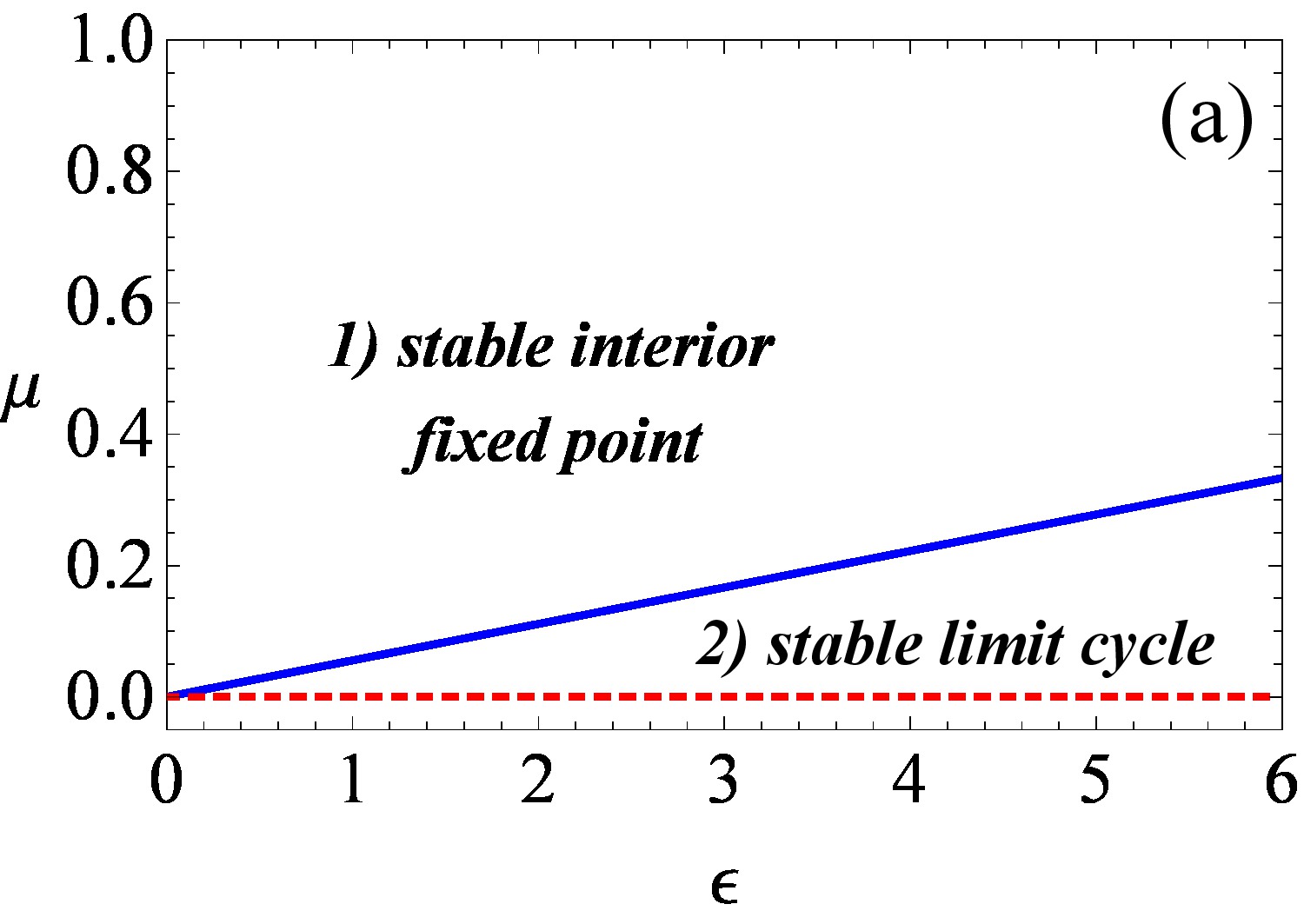}
\vspace{0.1cm}
\includegraphics[scale = 0.33]{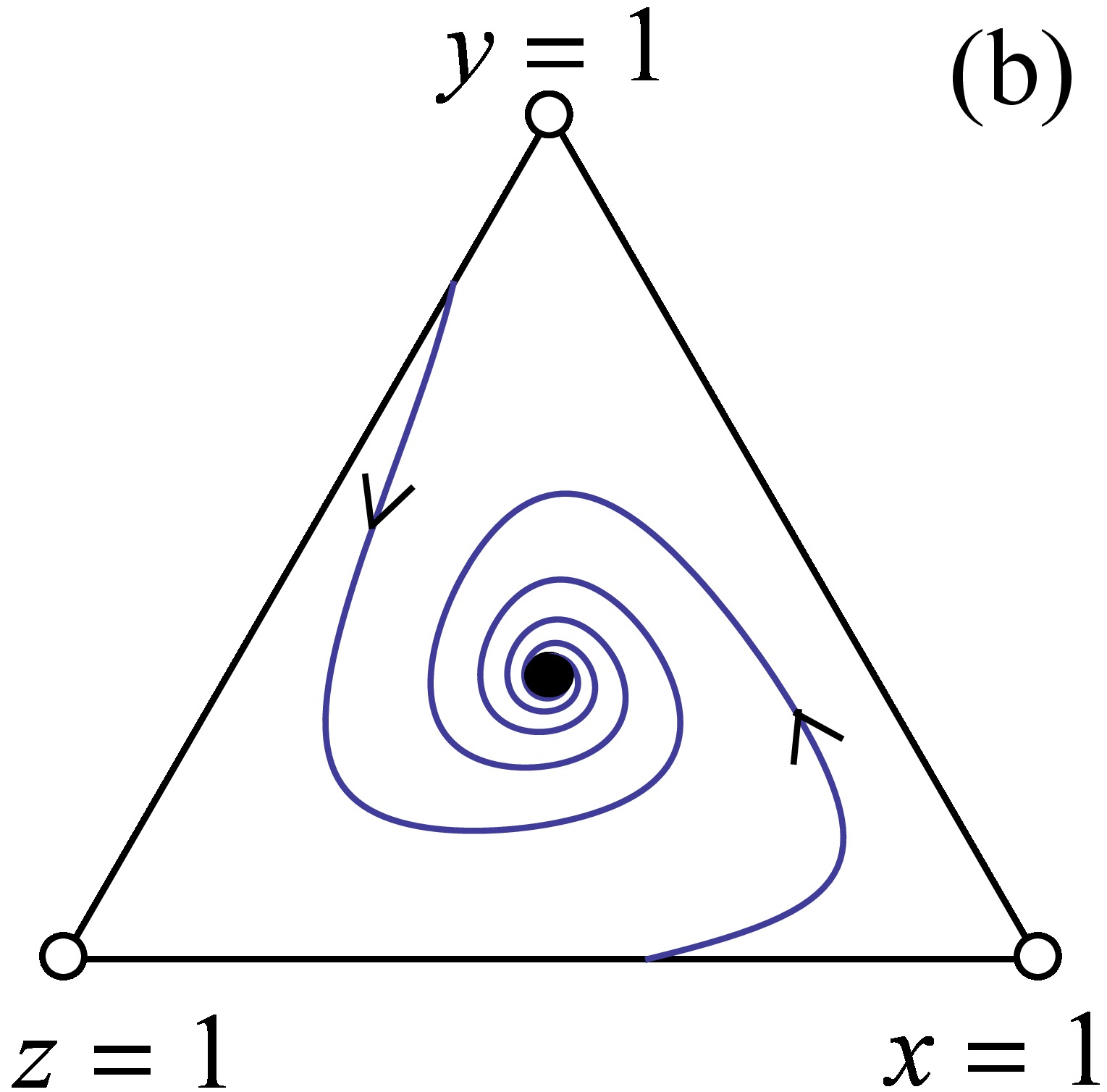} \includegraphics[scale = 0.33]{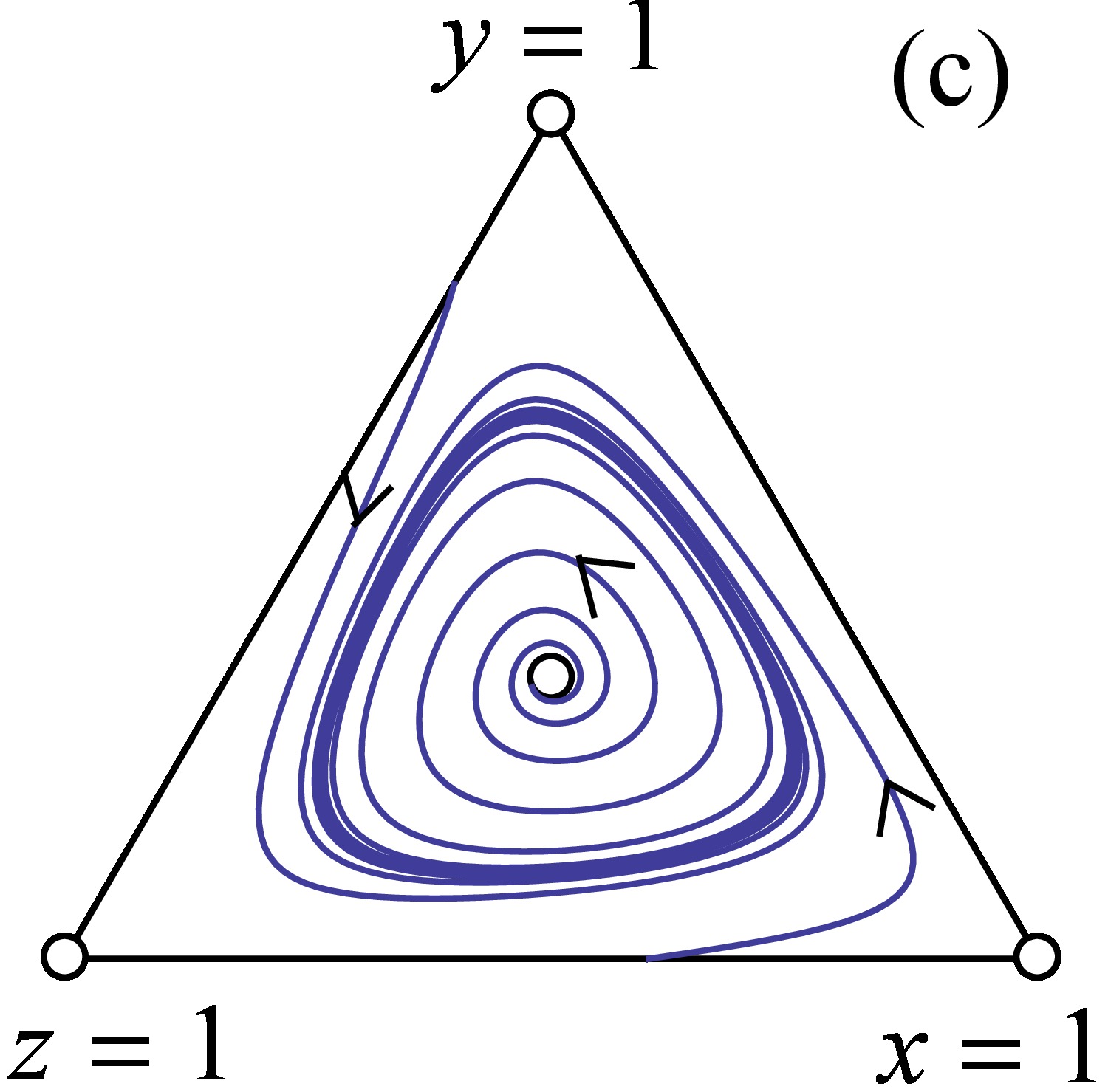}

\caption{\label{fig:RPSglobal} (Color online) Stability diagram for Eq.~\eqref{eqn:globalrighttriangle} and its associated phase portraits, in the case of  global mutation. (a) Stability diagram. Locus of supercritical Hopf bifurcation, solid  blue curve; saddle connection, dotted red line. (b) Phase portrait corresponding to Region 1 of the stability diagram. The system has a stable interior point, corresponding to coexistence of all three strategies. (c) Phase portrait corresponding to Region 2. The system has a stable limit cycle.}
\end{figure}

Equation~\eqref{eqn:globalrighttriangle} has four fixed points: $\left(x^*,y^*,z^*\right)= \left(0,0,0\right), \left(1,0,0\right), \left(0,1,0\right)$, and $\left(\frac{1}{3}, \frac{1}{3}, \frac{1}{3}\right)$. The inner fixed point undergoes a Hopf bifurcation when $\displaystyle \mu_h = \epsilon/18$, as shown in~\cite{mobilia2010oscillatory}. Moreover, the system undergoes three simultaneous saddle connections when $\mu = 0$, meaning that when $\mu = 0$, there is a heteroclinic cycle linking the three corners of the simplex. Figure~\ref{fig:RPSglobal} shows the stability diagram of Eq.~\eqref{eqn:globalrighttriangle} as well as phase portraits corresponding to the different regions of Fig.~\ref{fig:RPSglobal}(a).


Incidentally, the occurrence of the heteroclinic cycle is not an artifact of the global mutation structure assumed here. For any mutation pattern, $\mu = 0$ will always (trivially) yield a heteroclinic cycle linking the three corners of the simplex, since that cycle is known to be present in the absence of mutation~\cite{mobilia2010oscillatory, hofbauer1998evolutionary,szabo2004rock,szolnoki2014cyclic,nowak2006evolutionary}. Thus, all the subsequent stability diagrams in this paper will show a saddle connection curve along the line $\mu = 0$.


\section{Single Mutations}
Because of the cyclic symmetry of  the Rock-Paper-Scissors game, it suffices to consider two of the six possible single-mutation pathways. So without loss of generality, we restrict attention to rock ($x$) $\overset{\mu}{\rightarrow}$ paper ($y$) and paper $(y)$ $\overset{\mu}{\rightarrow}$ rock $(x)$. 

The two cases are qualitatively different. In the first case, the direction of mutation reinforces the system's inherent tendency to flow from rock to paper. (Recall that paper beats rock, so trajectories tend to flow from rock to paper, as the paper population grows at rock's expense.) By contrast, in the second case, the mutation pathway runs counter to this natural flow.  

\subsection{Single mutation: $x \xrightarrow{\mu} y$ }

\begin{figure}[h!]

\center
\includegraphics[scale=0.5]{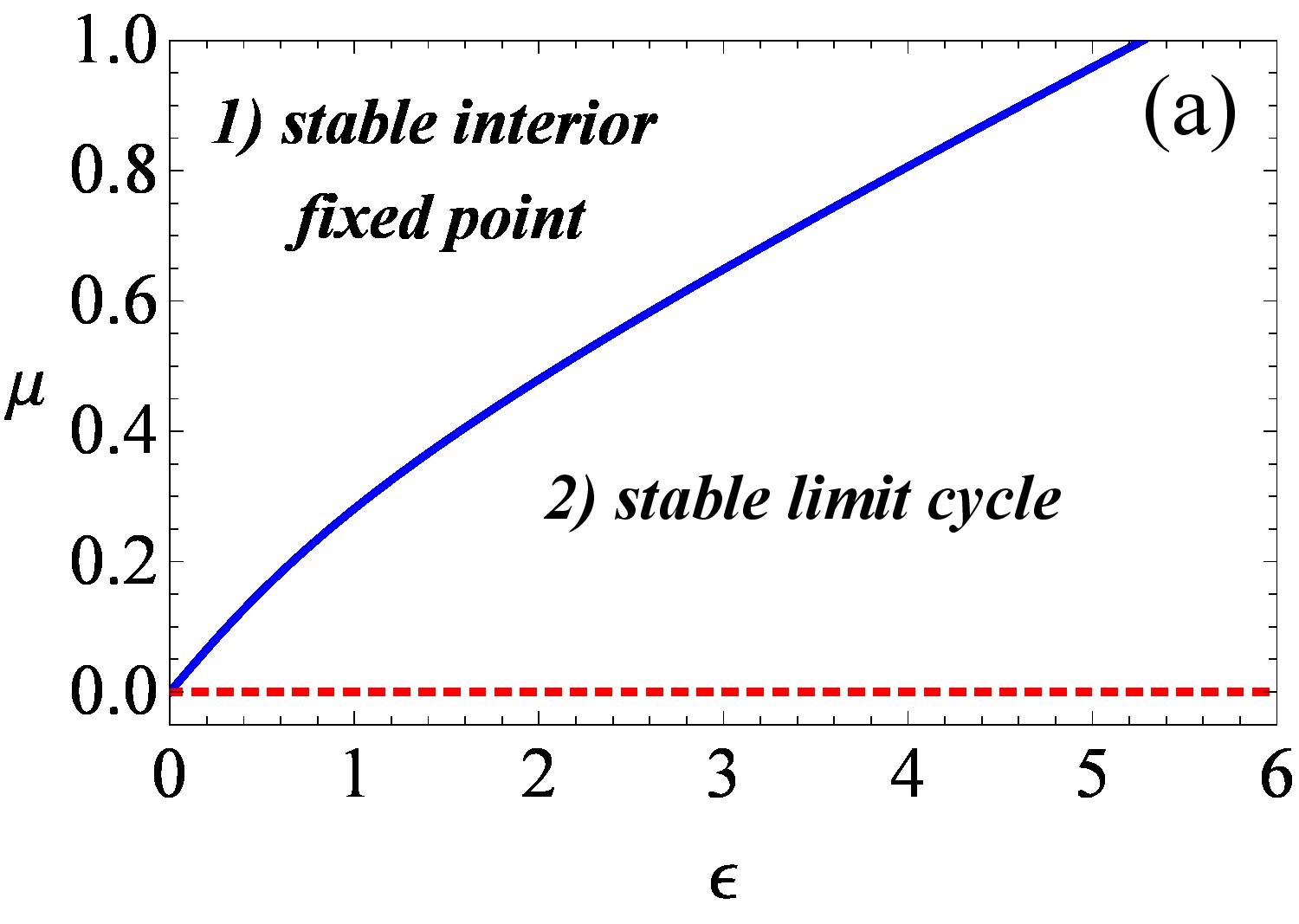}
\vspace{0.1cm}
\includegraphics[scale = 0.33]{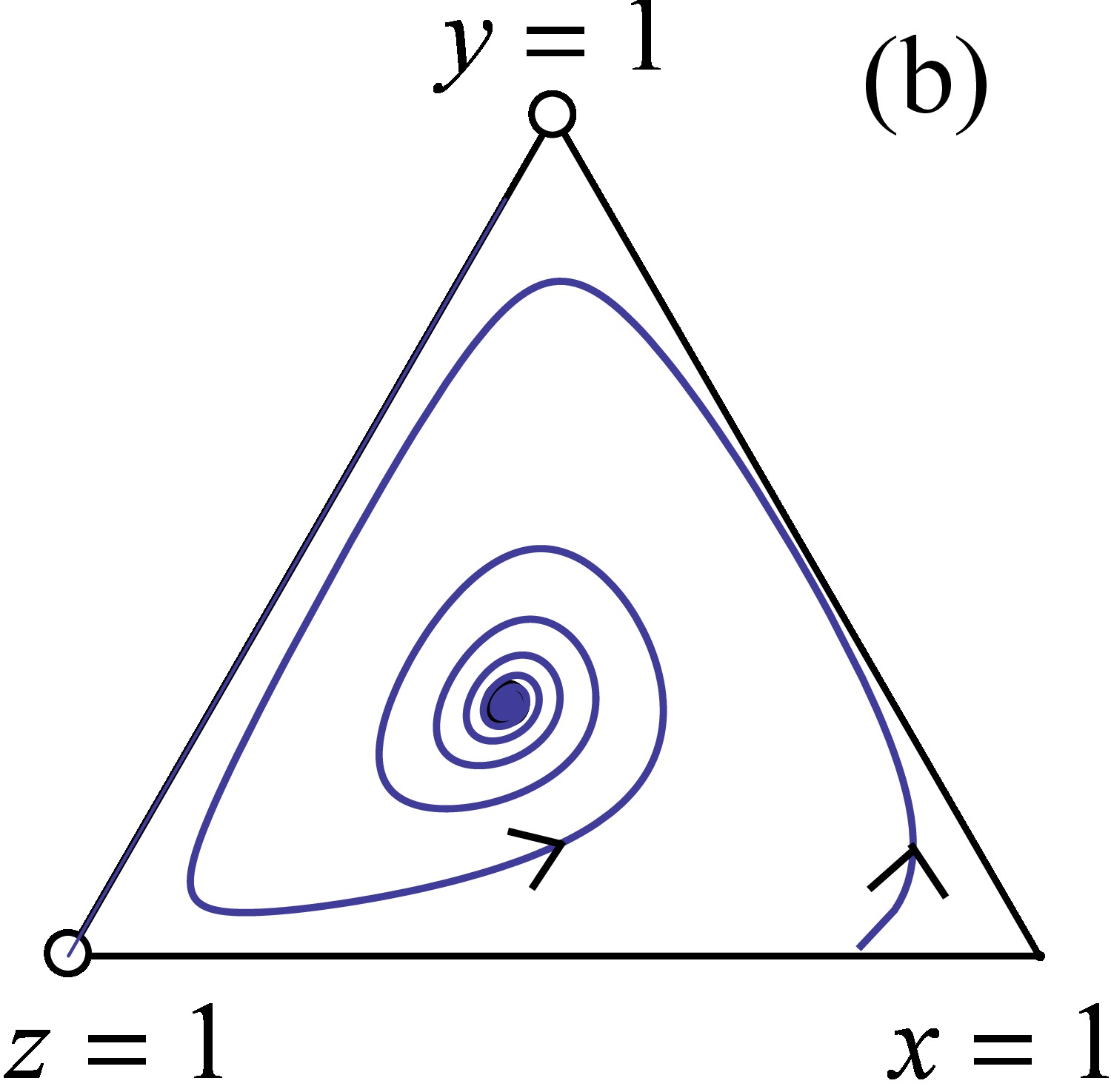} \includegraphics[scale = 0.33]{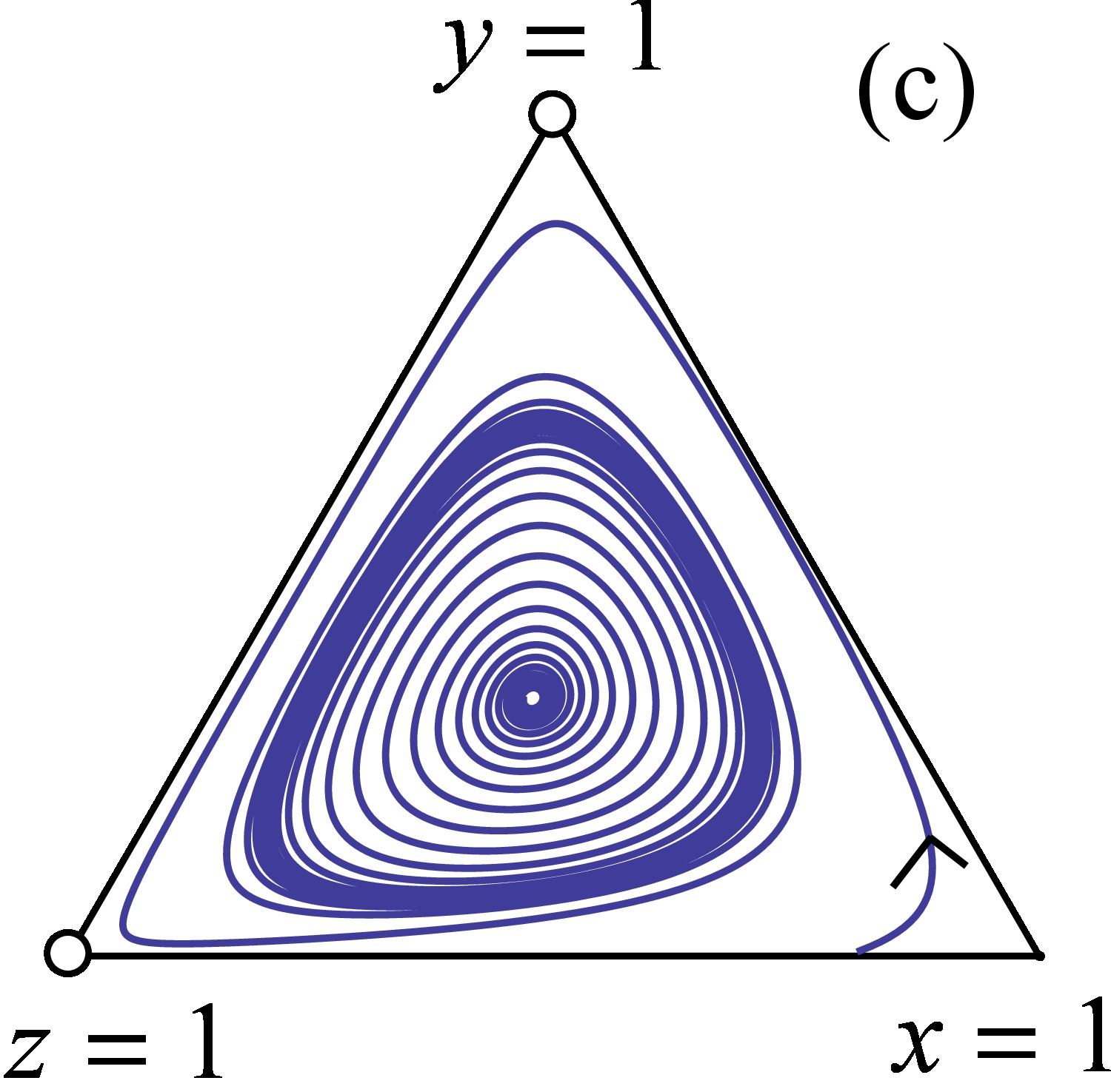}

\caption{\label{fig:RP} (Color online) Stability diagram and phase portraits for Eq.~\eqref{eqn:xtoy}. In this single-mutation case, rock $x$ mutates into paper $y$ at a rate $\mu$: $x\overset{\mu}{\rightarrow}y$ . (a) Stability diagram. Supercritical Hopf bifurcation, solid blue curve; saddle connection, dotted red line. (b) Phase portrait corresponding to Region 1 of the stability diagram. The system has a stable interior point, corresponding to coexistence of all three strategies. (c) Phase portrait corresponding to Region 2. The system has a stable limit cycle.
}
\end{figure}

When rock mutates into paper, the system becomes
\begin{eqnarray}
\dot{x} &=& x \left(f_x - \phi\right) -\mu x \nonumber\\
\dot{y} &=& y   \left(f_y - \phi\right) +\mu x.   \label{eqn:xtoy}
\end{eqnarray}
Figure~\ref{fig:RP} plots the stability diagram and phase portraits of Eq.~\eqref{eqn:xtoy}.


To derive the results shown in Fig.~\ref{fig:RP}, we observe first that Eq.~\eqref{eqn:xtoy} has three fixed points: $(0,0), (0,1),$ and an inner fixed point $(x^*_3,y^*_3)$, with
\begin{align*}
 x^*_3 &= \frac{\left(\epsilon +3\right) A_1+\epsilon  \left(3 \mu +\epsilon ^2+3 \mu  \epsilon -6\right)-9}{6 \left(\epsilon
    \left(\epsilon +3\right)+3\right)},\\
 y^*_3 &= \frac{-6 \mu +A_1+\epsilon  \left(-3 \mu +\epsilon +3\right)+3}{6 \left(\epsilon  \left(\epsilon +3\right)+3\right)},
\end{align*}
where the quantity $A_1$ is given by $$\displaystyle A_1=\sqrt{-3 \mu ^2 \epsilon ^2-6 \mu  \epsilon  (\epsilon  (\epsilon +3)+3)+(\epsilon  (\epsilon +3)+3)^2}.$$

These fixed points display some notable differences from those found above, when mutation was either absent or global. For example, the state in which only rock exists, corresponding to the corner $(x,y,z)=(1,0,0)$ of the simplex, is no longer a fixed point, since $x$ is now constantly mutating into $y$. So paper must exist whenever rock does. 
The complicated interior fixed point
$(x^*_3,y^*_3)$ can be regarded as a perturbation of $\left(\frac{1}{3}, \frac{1}{3}\right)$, in the sense that $(x^*_3,y^*_3) = (\frac{1}{3}+\delta_x, \frac{1}{3}+\delta_y)$ where both $\delta_x$ and $\delta_y$ approach 0 as $\mu\rightarrow 0$.

Equation~\eqref{eqn:xtoy} produces stable limit cycles when the interior fixed point undergoes a supercritical Hopf bifurcation. This transition can be shown to occur at
\begin{eqnarray}
\begin{split}
 \mu_h &=\frac{2 \left(\sqrt{\epsilon  (\epsilon +2) (4 \epsilon  (\epsilon +2)+9)+9}-3\right)-3 \epsilon  (\epsilon +2)}{7 \epsilon}\\
            & \approx  \frac{\epsilon }{3}-\frac{4 \epsilon ^3}{27}+\frac{4 \epsilon ^4}{27}-\frac{4 \epsilon ^5}{243}+O\left(\epsilon ^{6}\right).
\nonumber
\end{split}
\end{eqnarray}
The stable limit cycle created in the Hopf bifurcation grows into a heteroclinic cycle as $\mu \rightarrow 0$. 

\subsection{Single mutation: $y\xrightarrow{\mu}x$}
In the case where paper $y$ mutates into rock $x$ at a rate $\mu$, the system becomes
\begin{eqnarray}
\dot{x} &=& x (f_x - \phi) -\mu x \nonumber\\
\dot{y} &=& y   (f_y - \phi) +\mu x.   \label{eqn:ytox}
\end{eqnarray}

\begin{figure}[h!]

\center
\includegraphics[scale=0.5]{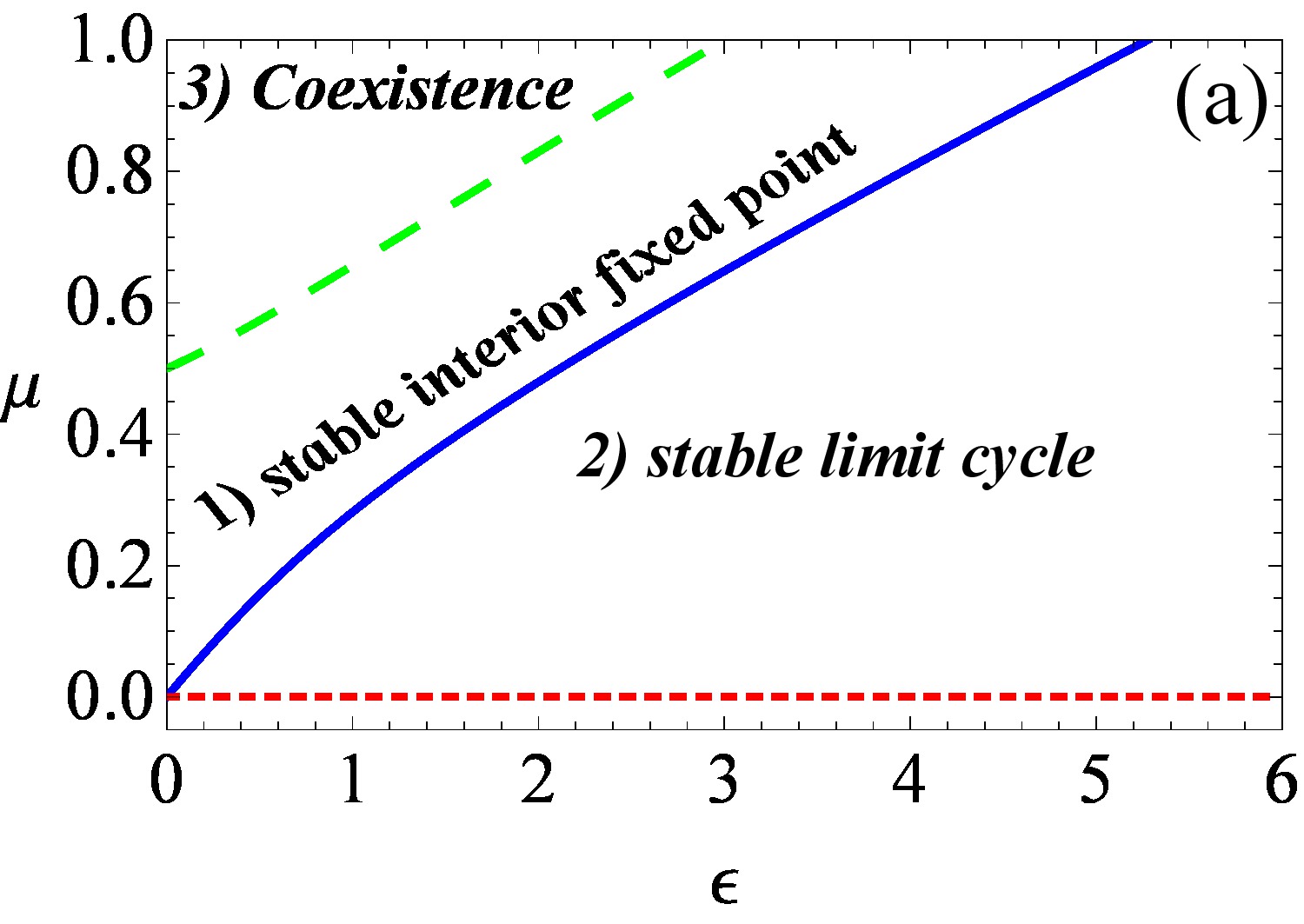}

\center
\includegraphics[scale = 0.33]{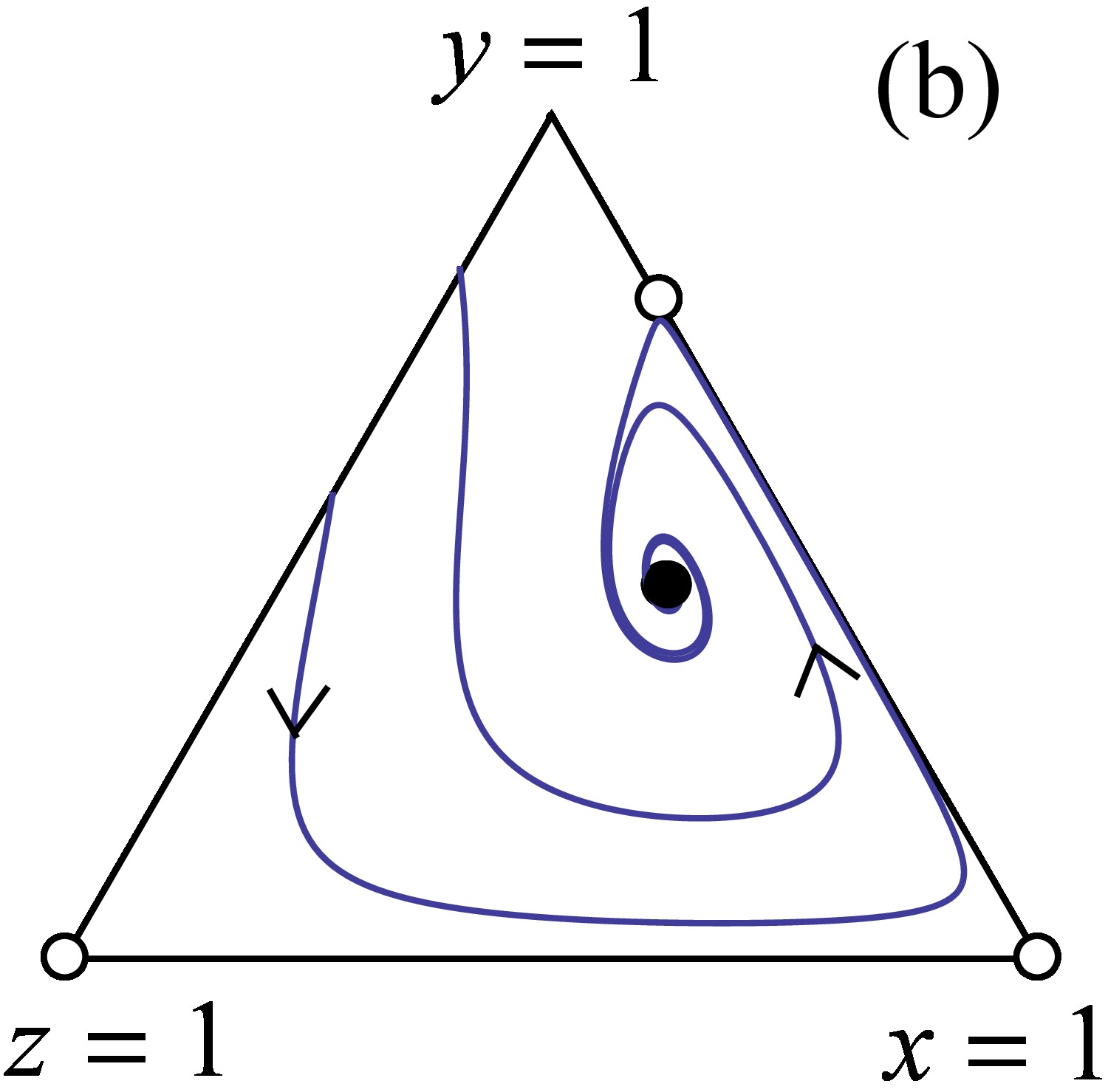} \includegraphics[scale = 0.33]{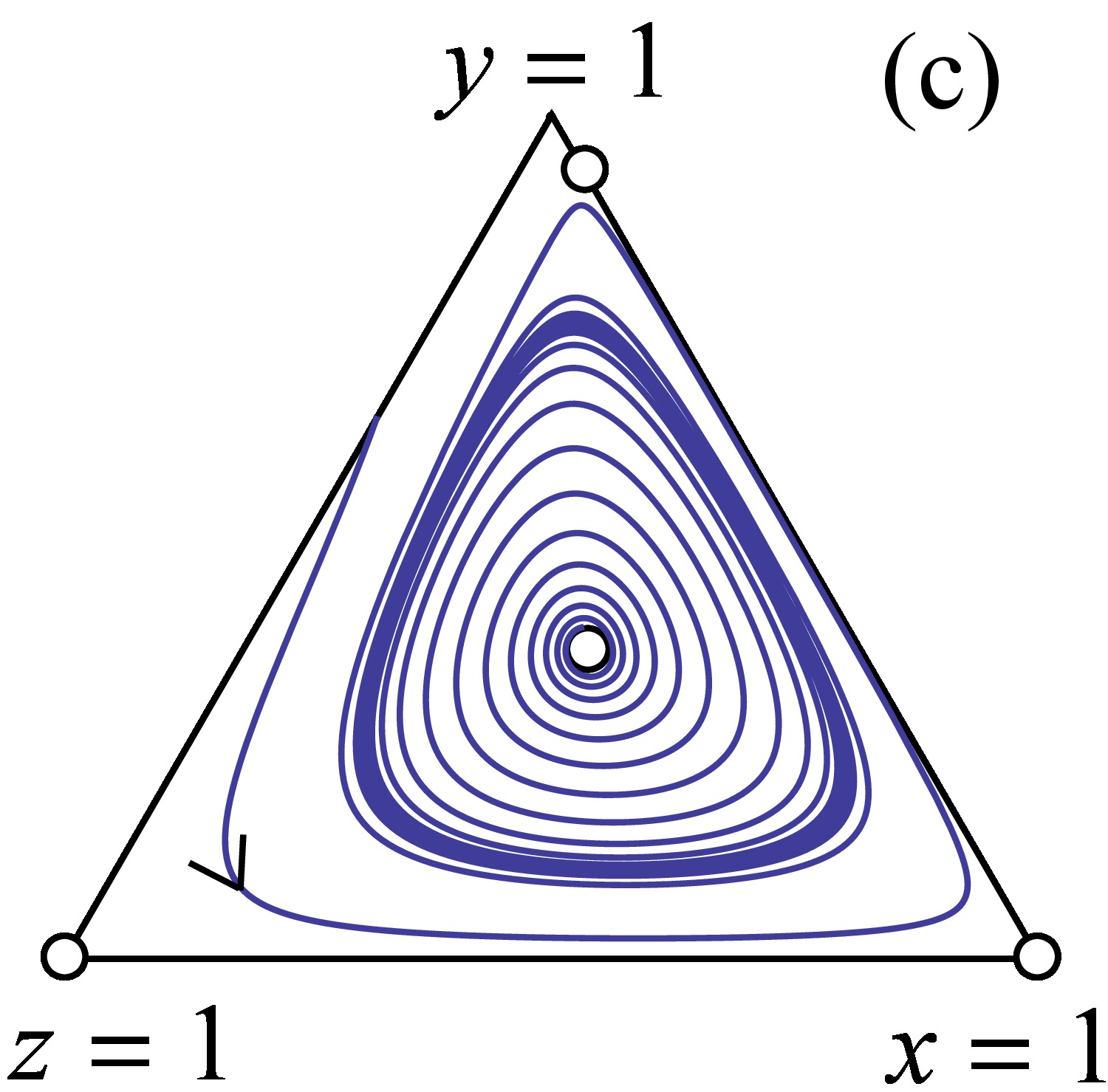}

\center
\includegraphics[scale = 0.33]{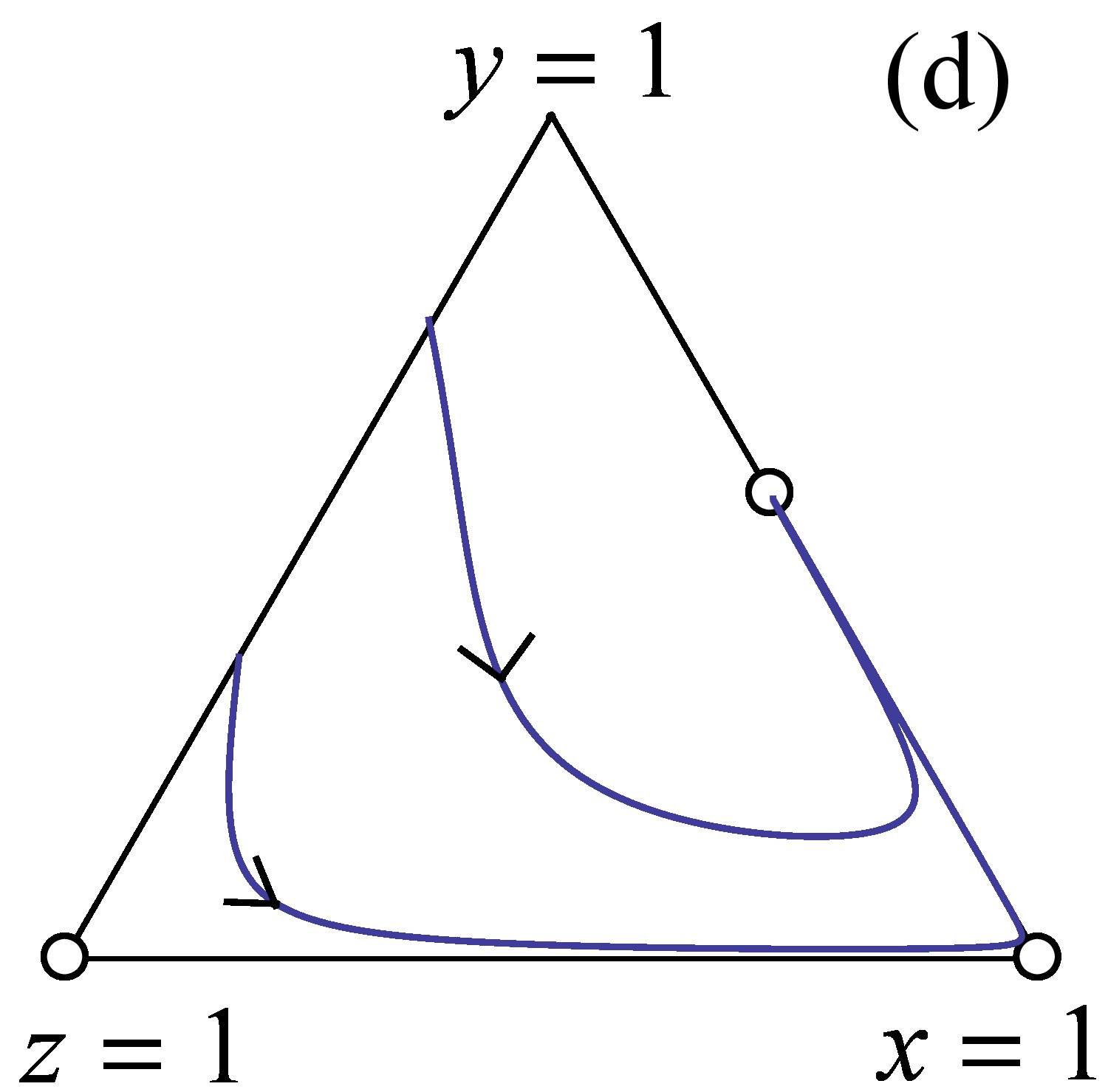}

\caption{\label{fig:PR} (Color online) Stability diagram and phase portraits for Eq.~\eqref{eqn:ytox}, corresponding to a single mutation pathway in which paper $y$ mutates into rock $x$ at a rate $\mu$. Note that this direction of mutation goes ``against the flow'' in the following sense: in the absence of mutation, paper beats rock and thus the flow normally tends to convert $x$ into $y$. Thus, the direction of mutation assumed here opposes that flow. (a) Stability diagram. Supercritical Hopf bifurcation, solid blue curve; saddle connection, dotted red line; transcritical bifurcation, dashed green curve. In the coexistence region above the transcritical bifurcation curve, rock coexists with paper at the stable fixed point. (b) Phase portrait corresponding to Region 1 of the stability diagram. The system has a stable interior point, corresponding to coexistence of all three strategies. (c) Phase portrait corresponding to Region 2. The system has a stable limit cycle. (d) Phase portrait corresponding to Region 3, where rock and paper coexist and scissors has gone extinct. The system has a stable point on the boundary line joining rock and paper.} 
\end{figure}

Figure~\ref{fig:PR} plots the stability diagram and phase portraits of Eq.~\eqref{eqn:ytox}. Note the existence of a new region in parameter space, bounded below by a transcritical bifurcation curve. This region did not exist when the mutation was in the direction of the system's inherent flow, depicted earlier in Fig.~\ref{fig:RP}. (The existence of such a region also holds if there are multiple mutations, unless the inner fixed point is $\left(\frac{1}{3},\frac{1}{3}\right)$, as we will see below.)


The results shown in Fig.~\ref{fig:PR} can be derived analytically, as follows. Because the direction of mutation $y\overset{\mu}{\rightarrow}x$ here opposes the system's inherent flow, the fixed points of Eq.~\eqref{eqn:ytox} now include a new fixed point $(x^*_4,y^*_4)$ on the boundary line given by $z=0$ and $x+y=1$.
The fixed points are thus $(0,0), (1,0),$ an interior fixed point $(x^*_3,y^*_3)$ similar to that found earlier, and a fourth fixed point on the boundary.  The coordinates of the nontrivial fixed points are given by
\begin{align*}
 x^*_3 &= \frac{6 \mu +A_1+\epsilon  (3 \mu +\epsilon +3)+3}{6 (\epsilon  (\epsilon +3)+3)}, \\
 y^*_3 &=\frac{\left(-2 \epsilon -3\right) A_1+\epsilon  (-3 \mu +\epsilon  (4 \epsilon +15)+21)+9}{6 \epsilon  (\epsilon  (\epsilon +3)+3)},\\
 x^*_4 &=\frac{\epsilon +1+\sqrt{(\epsilon +1)^2-4 \mu  \epsilon }}{2 \epsilon },\\
 y^*_4 &=\frac{\epsilon -1-\sqrt{(\epsilon +1)^2-4 \mu  \epsilon }}{2 \epsilon },
\end{align*}
where $A_1$ is given by the expression obtained earlier. 

As shown in Fig.~\ref{fig:PR}(a), the fixed point $(x^*_4,y^*_4)$ exists only for parameter values above the transcritical bifurcation curve.  By linearizing about the fixed point and seeking a zero-eigenvalue bifurcation, we find that the transcritical bifurcation curve is given by $$ \mu_{trans} =\frac{\epsilon -\sqrt{\epsilon }+1}{\sqrt{\epsilon }+1}.$$

Interestingly,  the Hopf bifurcation curve for Eq.~\eqref{eqn:ytox} is identical to that for Eq.~\eqref{eqn:xtoy}, even though the two systems have qualitatively different dynamics. We have no explanation for this coincidence. It does not follow from any obvious symmetry, as evidenced by the fact that the inner fixed points $(x^*_3,y^*_3)$ are different in the two cases.


\section{Double Mutations}
\begin{figure}[h!]

\center
\includegraphics[width=50mm,height=3.5cm,scale=0.4]{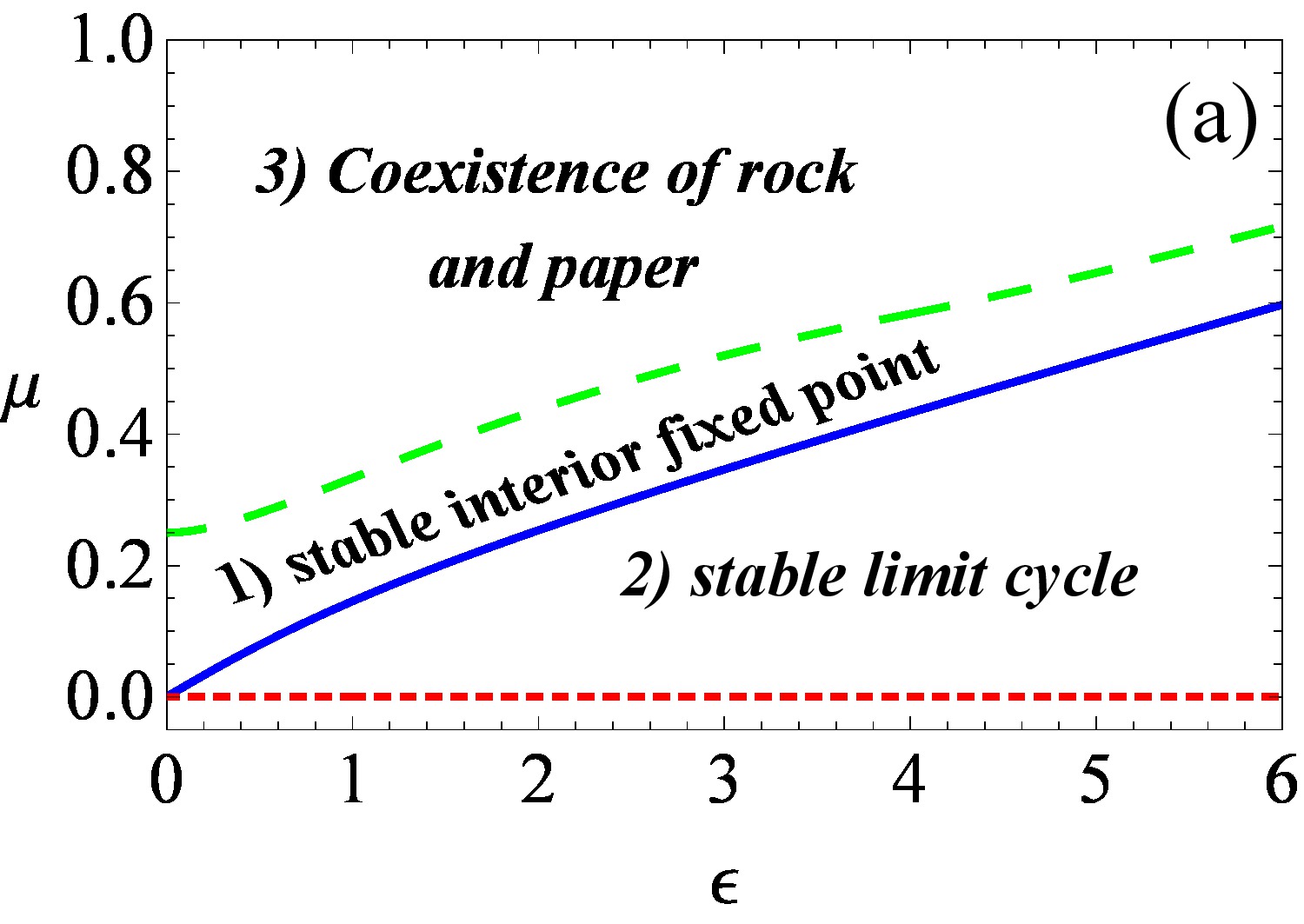}

\center
\includegraphics[width=50mm,height=3.5cm]{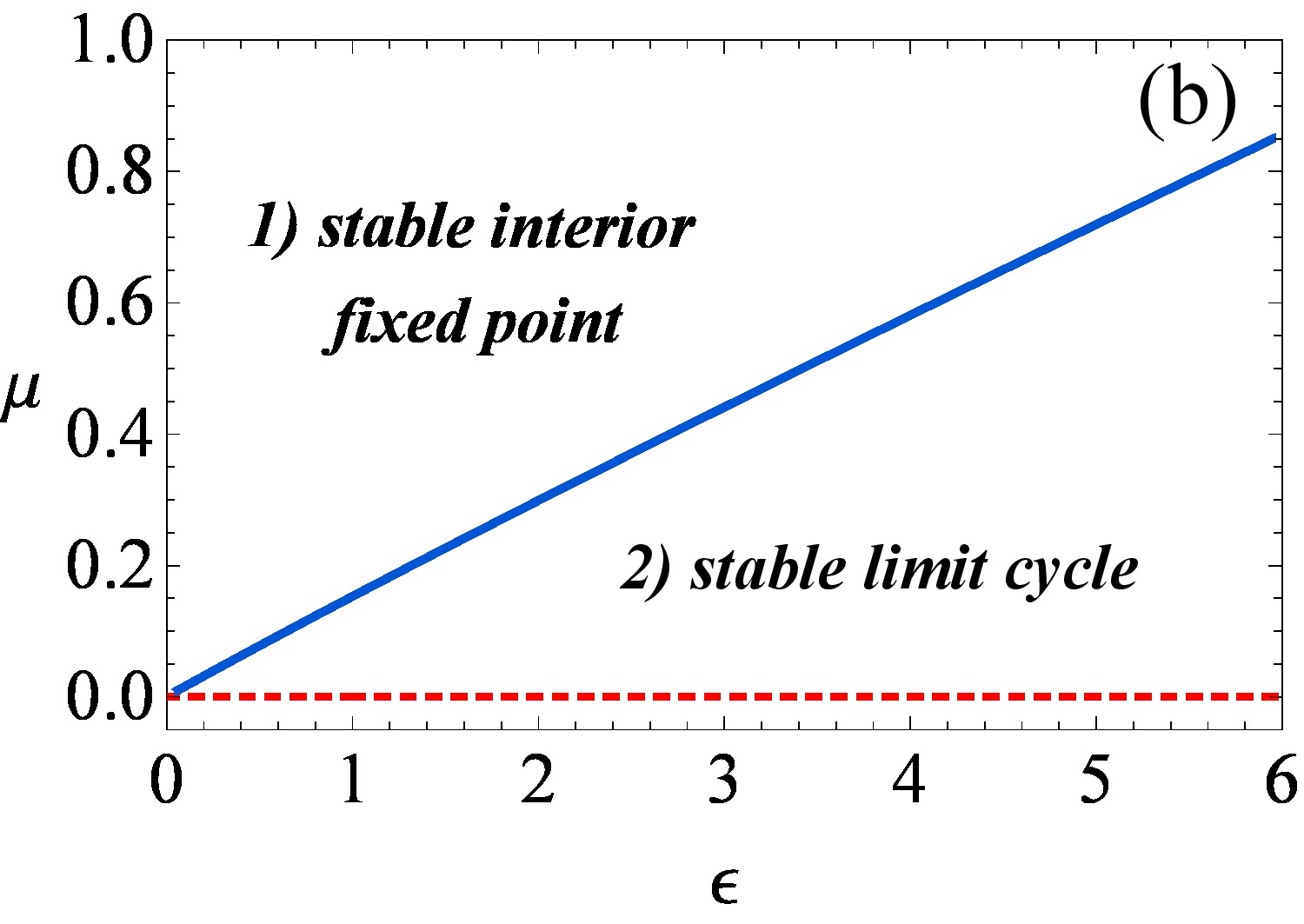}

\center
\includegraphics[width=50mm,height=3.5cm,scale=0.4]{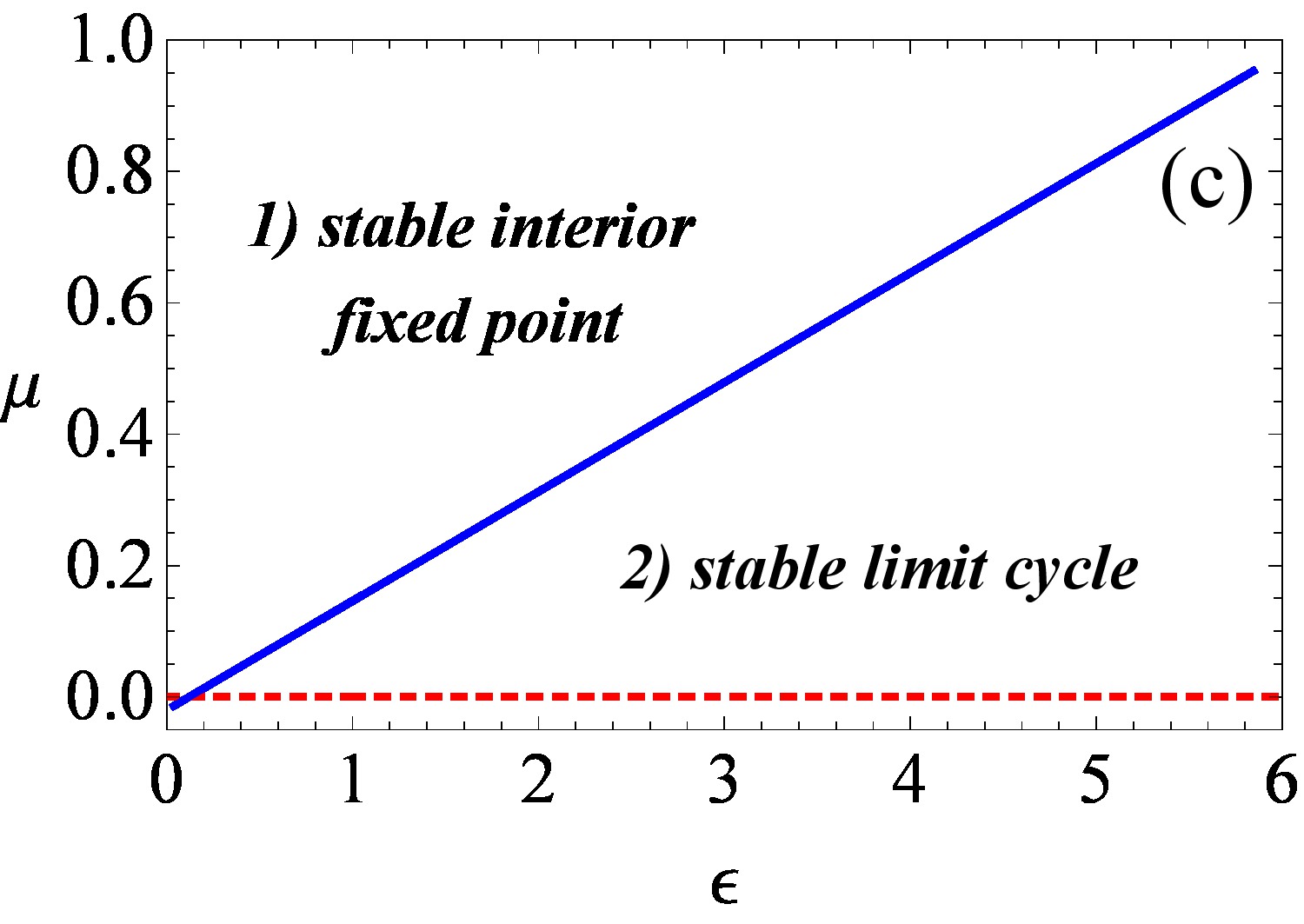}

\caption{\label{fig:double} (Color online) Stability diagram of the replicator equations for different patterns of two mutations, both occurring at a rate $\mu$. Supercritical Hopf bifurcation, solid blue curve;  transcritical bifurcation, dashed green curve; saddle connection, dotted red line. (a) Opposing mutations: $z\overset{\mu}{\rightarrow}y$ and $x\overset{\mu}{\rightarrow}y$. (b) Mutations in the same direction of circulation: $z\overset{\mu}{\rightarrow}x$ and $x\overset{\mu}{\rightarrow}y$. (c) Bidirectional mutation between the same two species: $x\overset{\mu}{\rightarrow}y$ and $y\overset{\mu}{\rightarrow}x$.}
\end{figure}

If we allow mutations to occur along two pathways instead of one, and assume that they both occur at the same rate $\mu$, then by the cyclic symmetry of the Rock-Paper-Scissors game there are four qualitatively different cases to consider. The analysis becomes more complicated than with single mutations, so we omit the details and summarize the main results in the stability diagrams shown in Figure~\ref{fig:double}. The key point is that in every case, stable limit cycles exist arbitrarily close to the origin $(\epsilon, \mu) = (0,0)$ in parameter space, consistent with the conjecture discussed in the Introduction.

Figure~\ref{fig:double}(a) shows the stability diagram for the first case, defined by having one of the two mutations go in the direction of the system's inherent flow and the other against it. An example is $z\overset{\mu}{\rightarrow}y$ and $x\overset{\mu}{\rightarrow}y$. As shown in Fig.~\ref{fig:double}(a), the resulting stability diagram has three regions and resembles what we saw earlier in Fig.~\ref{fig:PR}(a). Nothing qualitatively new happens if both pathways go against the flow so this case is omitted. 

The second case occurs when both mutation pathways go in the same direction relative to the flow, as in Fig.~\ref{fig:double}(b). Then the stability diagram has only two regions, similar to Fig.~\ref{fig:RP}(a). 

The final case, shown in Fig.~\ref{fig:double}(c), occurs when the two mutation pathways go in opposite directions between the same two species, as in $x\overset{\mu}{\rightarrow}y$ and $y\overset{\mu}{\rightarrow}x$. Again, the stability diagram shows only two regions. The boundary between the regions turns out to be exactly straight in cases like this. Specifically, we find that the Hopf curve here is given by the line $\mu_h = \epsilon/6$. The key to the analysis is the observation that the interior fixed point reduces to $x=y=z=\frac{1}{3}$ for this case. This convenient symmetry property eases the calculation of the Hopf bifurcation curve. Indeed, the Hopf curve continues to be straight, even for more complex patterns of mutation, as long as all three species are equally populated at the interior fixed point, as we will show next.


\begin{table}[h!]
\caption{Hopf curve for different mutation patterns}
\label{table:Hopf}
\begin{ruledtabular}
\begin{tabular}{ccccc}
   Mutation patterns\\ (rate $\mu$)& \# of mutations    &   Hopf curve
  \\ \hline\\
\includegraphics[trim = 0cm 18cm  15cm 6cm,clip, scale=0.4]{arrows.jpg}&2&  $\displaystyle \mu_h=\frac{\epsilon}{6}$
  \\ \hline\\
\includegraphics[trim = 0cm 14cm  15cm 10cm,clip, scale=0.4]{arrows.jpg} &3&  $\displaystyle \mu_h=\frac{\epsilon}{9}$
  \\ \hline\\
\includegraphics[trim = 0cm 6cm  15cm 18cm,clip, scale=0.4]{arrows.jpg}&4&  $\displaystyle \mu_h=\frac{\epsilon}{12}$
  \\ \hline\\

\includegraphics[trim = 0cm 10cm  15cm 14cm,clip, scale=0.4]{arrows.jpg} & 6&  $\displaystyle \mu_h=\frac{\epsilon}{18}$
  \\ 
\end{tabular}
\end{ruledtabular}
\end{table}


\section{Complex Mutations}
The analytical treatment of the model becomes prohibitively messy as one adds more mutation pathways. To make further progress, let us restrict attention to mutation patterns that preserve $(x,y) = (\frac{1}{3},\frac{1}{3})$ as the inner fixed point for all values of $\epsilon$ and $\mu$.  These mutation patterns are shown in Table~\ref{table:Hopf}. We computed the Hopf curve analytically for all of them and noticed something curious: The equation of the Hopf curve is always $$ \mu_h=\frac{\epsilon}{3 \times (\# \text{ of mutations)}},$$ as shown in Table~\ref{table:Hopf}.

To derive this formula, we write the replicator equations for all the mutation patterns in Table~\ref{table:Hopf} as a single system:
\begin{eqnarray}
\dot{x} &=& x (f_x - \phi) + \mu(-[\alpha_y x+\alpha_z x]+\beta_x y+\gamma_x z) \nonumber\\
\dot{y} &=& y   (f_y - \phi) + \mu(-[\beta_x y+\beta_z y]+\alpha_y x+\gamma_y z) \nonumber\\
\dot{z} &=& z   (f_z - \phi) + \mu(-[\gamma_x z+\gamma_y z]+\alpha_z x+\beta_z y)   \label{eqn:all}.
\end{eqnarray}
Here the indicator coefficients given by the various $\alpha, \beta,$ and $\gamma$ are set to 0 or 1, depending on which mutation pathways are absent or present. For example, we set $\alpha_y = 1$ if $x$ mutates into $y$. Otherwise, we set $\alpha_y = 0$. The same sort of reasoning applies to the various $\beta_j$ and $\gamma_k$. 

To ensure that $(\frac{1}{3},\frac{1}{3})$ is a fixed point of Eq.~\eqref{eqn:all}, the indicator coefficients must satisfy certain algebraic constraints. These are given by
\begin{eqnarray}
 \beta_x + \gamma_x = \alpha_y + \alpha_z\nonumber\\
\alpha_y + \gamma_y= \beta_x+\beta_z \nonumber\\ 
\alpha_z + \beta_z  =  \gamma_x + \gamma_y \label{eqn:abc}.
\end{eqnarray}
Using Eqs.~\eqref{eqn:all} and ~\eqref{eqn:abc},  one can then show that the fixed point $(\frac{1}{3},\frac{1}{3})$ undergoes a Hopf bifurcation at $$ \mu_h=\frac{\epsilon}{3(\alpha_y+\alpha_z+\beta_x+\beta_z+\gamma_x+\gamma_y)},$$ which yields the results in Table~\ref{table:Hopf}.

The upshot is that a curve of supercritical Hopf bifurcations emanates from the origin in parameter space. Hence, for complex mutation patterns that preserve the fixed point $(x,y) = (\frac{1}{3},\frac{1}{3})$, we have confirmed the conjecture that stable limit cycles exist for parameters arbitrarily close to the zero-sum, zero-mutation-rate limit of the replicator equations for the Rock-Paper-Scissors game.


\section{Discussion}
Our main result is that for a wide class of mutation patterns, the replicator-mutator equations for the Rock-Paper-Scissors game have stable limit cycle solutions. For this class of mutation patterns, a tiny rate of mutation and a tiny departure from a zero-sum game is enough to destabilize the coexistence state of a Rock-Paper-Scissors game and to set it into self-sustained oscillations. 

However, we have not proven that limit cycles exist for \emph{all} patterns of mutation. That question remains open.  

Another caveat is that our results have been obtained for one version of the replicator-mutator equations, namely, that in which the mutation terms are \emph{added} to the replicator vector field. In making this choice we are following Mobilia~\cite{mobilia2010oscillatory}, who investigated the effect of additive global mutation on the replicator dynamics for the Rock-Paper-Scissors game. But there is another way to include the effect of mutation in the replicator equations: one can include it \emph{multiplicatively} as in, e.g., Ref.~\cite{nowak2006evolutionary}. In biological terms, the multiplicative case makes sense if the mutations occur in the offspring, whereas in the additive case they occur in the adults. We found the additive case easier to work with mathematically, but it would be interesting to see if our results would still hold (or not) in the multiplicative case.



%


\end{document}